\documentclass[runningheads]{llncs}
\usepackage{graphicx}
\usepackage{amssymb}
\usepackage{amsmath, mdframed}
\usepackage{comment}
\usepackage{enumerate}

\newcommand{\0}{\mathcal{O}}
\renewcommand{\a}{\alpha}
\renewcommand{\b}{\beta}
\newcommand{\e}{\epsilon}
\newcommand{\ff}{\mathsf{f}}
\newcommand{\F}{\mathbb{F}}
\renewcommand{\H}{\mathcal{H}}

\newcommand{\m}{\mathfrak{m}}
\renewcommand{\P}{\mathbb{P}}
\newcommand{\rk}{\textnormal{rk}}
\newcommand{\se}{\subseteq}
\renewcommand{\v}{\mathrm{v}}
\newcommand{\Z}{\mathbb{Z}}

\begin{document}
\title{The group structure of elliptic curves over $\Z/N\Z$}
%
%
\author{Massimiliano Sala\inst{1}\orcidID{0000-0002-7266-5146} \and
Daniele Taufer\thanks{Supported in part by the European Union’s H2020 Programme under grant agreement number ERC-669891, and in part by the Research Foundation - Flanders (FWO), project 12ZZC23N and travel grant V425623N.}\inst{2}\orcidID{0000-0003-3402-4863}}
\authorrunning{M. Sala et D. Taufer}
%
\institute{University of Trento, Italy. \email{massimiliano.sala@unitn.it} \and
KU Leuven, Belgium.
\email{daniele.taufer@kuleuven.be}}
\maketitle              
\begin{abstract}
We characterize the possible groups $E(\Z/N\Z)$ arising from elliptic curves over $\Z/N\Z$ in terms of the groups  $E(\F_p)$, with $p$ varying among the prime divisors of $N$.
This classification is achieved by showing that the infinity part of any elliptic curves over $\Z/p^e\Z$ is a $\Z/p^e\Z$-torsor, of which a generator is exhibited.
As a first consequence, when $E(\Z/N\Z)$ is a $p$-group, we provide an explicit and sharp bound on its rank.
As a second consequence, when $N = p^e$ is a prime power and the projected curve $E(\F_p)$ has trace one, we provide an isomorphism attack to the ECDLP, which works only by means of finite rings arithmetic.

\end{abstract}
\section{Introduction}
\label{S:1}
Elliptic curves have been providing number theory with a fertile field of intense research for the last century, from theoretic \cite{Mordell,Merel,Wiles,Modularity}, algorithmic \cite{Schoof,BosmaFact} and applied \cite{Miller,Koblitz,ECDSA,Shparlinski} sides.
In their basic definition, these objects consist of non-singular plane projective cubics, defined as the zero-set of a Weierstrass polynomial over a given base field.
It is well-known that these curves are actually abelian varieties with the chord-tangent sum
\cite{Silverman,Husemoller,Washington}.
The study of the group structure arising from this operation has attracted huge attention and its grasp has proven to be remarkably challenging.
Beyond its indisputable algebraic interest, the security of cryptographic protocols based on these curves relies upon the nature of their addition operation, hence the investigation of these groups has received impetus in the last decades.

When the underlying field is finite, any group that may be realized as the point group of an elliptic curve is known \cite{Ruck,Voloch}. Nevertheless, both their distribution
\cite{BPS} and their efficient explicit description \cite{GrpGenerators} are lines of open research.
We refer to \cite{GrpSurvey} for an overview of the known classification of groups arising from curves with a Weierstrass model.

Elliptic curves may also be defined over rings, among which $\Z/N\Z$ is a significant instance both from a theoretical perspective \cite{ECNTA} and for cryptographic applications \cite{LenstraFact,CryptoZNZ}.
In this paper, we are mainly interested in their algebraic, especially groupal, properties: we classify all the possible groups arising from elliptic curves over any residue ring $\Z/N\Z$ in terms of their projected components modulo the prime divisors of $N$.
More precisely, if $p$ is a prime integer and $\v_p(N)$ is the $p$-adic valuation of $N$, the Chinese Reminder Theorem provides a group isomorphism
\[
    E(\Z/N\Z) \simeq \bigoplus_{p | N} E(\Z/p^{\v_p(N)}\Z),
\]
whose components are known \cite{ECNTA,Washington} to split as
\[
    E(\Z/p^{\v_p(N)}\Z) = H \oplus E(\F_p).
\]
The subgroup at infinity $H$, given by the kernel of the canonical projection, is known to be a $p$-group, since $|H| = p^{\v_p(N)-1}$ \cite{ECNTA}.
However, the structure of this group was only recently determined in terms of $0$-layers of elliptic loops \cite{ELoops}.

In this work, we provide a complete classification result based only on the arithmetic of curves over $\Z/N\Z$, especially when projecting to anomalous elliptic curves.
In particular, we prove the following group isomorphism:
\[
    E(\Z/N\Z) \simeq \bigoplus_{\substack{p | N\\ |E(\F_p)| \neq p}} E(\F_p) \oplus \Z/p^{\v_p(N)-1}\Z \oplus \bigoplus_{\substack{p | N\\ |E(\F_p)| = p}} G_p,
\]
where every $G_p$ may be either $\Z/p^{\v_p(N)}\Z$ or $\F_p \oplus \Z/p^{\v_p(N)-1}\Z$.
This result is obtained by proving that the infinity part of $E(\Z/p^e\Z)$ is a $\Z/p^e\Z$-torsor, which is far from holding over generic local rings \cite{InvTau}.
By proving it, we refine the case $t=0$ of \cite[Proposition 10.3]{ELoops}, as we effectively exhibit the generator of this cyclic subgroup.

\medskip

From the above classification, we derive some consequences.
First, we give an explicit bound on the rank of $E(\Z/N\Z)$ when the points of such curve form a $p$-group. This bound is sharp and depends only on $p$, determining as a corollary infinitely many groups that cannot arise from such curves.
The proof of this bound also provides a systematic way for generating such $p$-curves of admissible ranks.
Second, we exhibit a polynomial-time isomorphism attack to the elliptic curve discrete logarithm problem (ECDLP) over anomalous curves.
Although similar attacks have already appeared \cite{SatohAraki,Smart}, we find this approach noteworthy as its correctness and execution may be elaborated with only finite rings arithmetic, which makes it slightly more elementary.


This paper is organized as follows. In Section \ref{S:2} we recall some known results and definitions, including the group structure of elliptic curves over finite fields and the definition of such curves over rings.
In Section \ref{S:3} the group structure of elliptic curves over $\Z/N\Z$ is investigated and we derive our main result (Theorem \ref{thm:classification}).
Consequently, in Section \ref{S:4} we present a bound to the rank of $p$-groups that may arise from elliptic curves over $\Z/N\Z$.
An isomorphism attack to the ECDLP over anomalous curves is described in Section \ref{S:5}.
Finally, conclusions and further work are discussed in Section \ref{S:6}.


\section{Preliminaries}
\label{S:2}


In this paper, $R$ always denotes a commutative ring with unity and $R^*$ is the set of its invertible elements.
We employ capital letters $X,Y,Z$ to denote elements of $R$, while lowercase ones are variables in $R[x,y,z]$.

\begin{definition}[Primitivity]
 A finite collection $\{X_i\}_{i \in \{0,\dots,n\}} \se R^{n+1}$ is called \emph{primitive} if the ideal $\langle \{X_i\}_{i \in \{0,\dots,n\}} \rangle_R$ is $R$ itself.
\end{definition}

\subsection{Elliptic curves over finite fields}

The trace $t$ of any elliptic curve over a finite field $\F_q$ is constrained by the Hasse bound \cite[Theorem V.1.1]{Silverman}, i.e.
\[
    t = q+1-|E(\F_q)| 
\]
is bounded by
\[
    -2\sqrt{q} \leq t \leq 2\sqrt{q}.
\]
Not every possible integer $t$ in the above interval occurs as the trace of an elliptic curve over $\F_q$, as detailed in \cite[Theorem 4.1]{Waterhouse}.
However, the same theorem shows that every such $t$ may be achieved if $q$ is a pure prime, i.e. the Hasse interval over prime fields is full.
From this work, a complete characterization of the possible point groups for elliptic curves over finite fields has seen the light, independently discovered by two authors \cite{Ruck,Voloch}.

    
    

By virtue of these works, we know all the possible groups arising from elliptic curves over finite fields, which we will use in Section \ref{S:3} to characterize those of curves over $\Z/N\Z$.

\subsection{Strong rank}

To deal with matrices over commutative rings, it is worth introducing a stronger notion of matrix rank.

\begin{definition}[Minor ideal]
Let $n,m \in \Z_{\geq 1}$ and $A \in M_{n,m}(R)$.
For every integer $1 \leq t \leq \min\{n,m\}$ we define the \emph{$t$-minor ideal} $I_t(A)$ as the ideal generated by the $t \times t$ minors of $A$.
We also define by convention $I_0(A) = R$ and for every $t > \min\{n,m\}$ we set $I_{t}(A) = (0)$.
\end{definition}
 
\begin{definition}[Strong rank]
Let $n,m \in \Z_{\geq 1}$ and $A \in M_{n,m}(R)$.
We define the \emph{strong rank} of $A$ as
\[
    \rk(A) = \max\{ t \in \Z_{\geq 0} \ | \ I_t(A) \neq (0) \}.
\]
\end{definition}

This notion of rank is easily shown to be never lower than the usual notion of rank over rings \cite[Chapter 4]{Brown}.
The convenience of using this rank relies on the following result.

\begin{lemma} \label{lem:rk}
Let $n,m \in \Z_{\geq 1}$ and $A \in M_{n,m}(R)$ be a matrix whose entries are primitive, then the following are equivalent.
\begin{enumerate}[$(i)$]
\item \label{rk:i} $\rk(A) = 1$.
\item \label{rk:ii} The $2 \times 2$ minors of $A$ vanish.
\item \label{rk:iii} All the primitive vectors of $R^n$ that may be obtained from an $R$-linear combination among the columns of $A$ are equal up to $R^*$-multiples. 
\end{enumerate}
\end{lemma}

\proof Let $A = (a_{i,k})_{\substack{1 \leq i \leq n \\ 1 \leq k \leq m}}$.\\
$[\ref{rk:i} \Rightarrow \ref{rk:ii}]$ Since $\rk(A) = 1$ then $I_2(A) = (0)$, hence all the generators of $I_2(A)$ vanish. \\
$[\ref{rk:ii} \Rightarrow \ref{rk:iii}]$ Let $v_1 = (v_{11}, \dots, v_{1n})$ and $v_2 = (v_{21}, \dots, v_{2n})$ be two primitive columns combinations. Since $v_1$ is primitive there are $\a_1, \dots, \a_n \in R$ with
\[
    \sum_{i = 1}^n \a_i v_{1i} = 1 \in R.
\]
Any $2 \times 2$ minor of the ($n \times 2$)-matrix $( v_1 | v_2 )$, whose columns are $v_1$ and $v_2$, is an $R$-linear combination of the $2 \times 2$ minors of $A$, hence it vanishes.
Thus, for every $i,j \in \{1, \dots, n\}$ we have $v_{1i}v_{2j} = v_{1j}v_{2i}$, then
\[
    v_2 = 1 \cdot v_2 = \left( \sum_{i = 1}^n \a_i v_{1i}v_{2j} \right)_{1 \leq j \leq n} = \left( \sum_{i = 1}^n \a_i v_{1j}v_{2i} \right)_{1 \leq j \leq n}
    = \left( \sum_{i = 1}^n \a_i v_{2i} \right) v_1.
\]
This proves that $v_2$ is a multiple of $v_1$, and since also $v_2$ is primitive then the scalar factor has to be a unit, i.e. $\sum_{i = 1}^n \a_i v_{2i} \in R^*$.\\
$[\ref{rk:iii} \Rightarrow \ref{rk:i}]$ For every pair of columns $c_k$ and $c_h$ of $A$ there is $r_{kh} \in R^*$ such that $c_h = r_{kh}c_k$. Therefore for every $1 \leq i,j \leq n$ we have
\[
    a_{ik}a_{jh} - a_{ih}a_{jk} = r_{kh} (a_{ik}a_{jk} - a_{ik}a_{jk}) = 0,
\]
which shows that $I_2(A) = (0)$.
Moreover, since the entries of $A$ are primitive we have $I_1(A) = R$, so that $\rk(A) = 1$.
\endproof

\subsection{Elliptic curves over rings}
\label{subsec:ECoverR}

Let $n$ be a non-negative integer.
The projective $n$-space over $R$ is defined in order to respect projections on any non-zero quotient of $R$, as follows.

\begin{definition}[Projective $n$-space]
 The \emph{projective $n$-space} over $R$ is the set of orbits of primitive tuples in $R^{n+1}$ under the action of elements $u \in R^*$ given by
 
 \[
    u (X_0, \dots, X_n) = (u X_0, \dots, u X_n).
\]
 It is denoted by $\P^n(R)$, while $(X_0 : \dots : X_n) \in \P^n(R)$ represents the orbit of $(X_0, \dots, X_n) \in R^{n+1}$.
\end{definition}

An elliptic curve over $R$ may be defined \cite{ECNTA} to properly extend a family of elliptic curves over $R/\m$, for $\m$ ranging among all the maximal ideals of $R$, provided that this ring satisfies the following condition.

\textbf{Condition} \cite{ECNTA}.
For every pair $n,m \in \Z_{\geq 1}$ and every matrix
\begin{equation} \label{R-Condition}
    A = (a_{ij})_{\substack{1 \leq i \leq n \\ 1 \leq j \leq m}} \in M_{n,m}(R)
\end{equation}
with strong rank $\rk(A) = 1$ and primitive entries, there exists an $R$-linear combination of the columns of $A$ whose entries are primitive.

In this work, we will only deal with elliptic curves that may be defined via their short Weierstrass equation, which is not restrictive when $6 \in R^*$.

\begin{definition}[Elliptic curve over $R$] Let $R$ be a commutative ring with unity satisfying Condition \ref{R-Condition} and let $A,B \in R$ such that 
\[
    \Delta_{A,B} = -(4A^3+27B^2) \in R^*.
\]
The \emph{elliptic curve} $E_{A,B}(R)$ is defined as
\[
    E_{A,B}(R) = \{(X:Y:Z) \in \P^2(R) \ | \ Y^2Z = X^3 + AXZ^2 + BZ^3 \}.
\]
Given an elliptic curve $E = E_{A,B}(R)$, we denote by $\0 = (0:1:0) \in E$ its \emph{zero element}, with $E^a = E \cap \P_{\text{aff}}^2(R)$ its \emph{affine points} and with $E^{\infty}$ the remaining points, which are called \emph{points at infinity}.
\end{definition}

On these curves a sum operation may be explicitly defined on an open covering of $E_{A,B}(R) \times E_{A,B}(R)$ by means of $(2,2)$-bidegree polynomials \cite{LanRup85,BosmaLenstra}.
This operation extends the usual point addition with respect to projections, i.e. for every proper ideal $I \subsetneq R$ we have a well-defined group homomorphism
\[
    \pi : E_{A,B}(R) \twoheadrightarrow E_{A,B}(R/I).
\]

We recall for convenience the two addition laws we employ in this work: the sum of $P_1 = (X_1:Y_1:Z_1)$ and $P_2 = (X_2:Y_2:Z_2)$ is given by any primitive linear combination of $(S_1 : S_2 : S_3)$ and $(T_1 : T_2 : T_3)$, where \footnote{Addition laws corresponding to $(0:0:1)$ and $(0:1:0)$ as in [\cite{BosmaLenstra}, Theorem 2]. }
\begin{align*} 
    S_1    = \ & (X_1Y_2-X_2Y_1)(Y_1Z_2 + Y_2Z_1) + (X_1Z_2-X_2Z_1)Y_1Y_2 \\
            & - A(X_1Z_2-X_2Z_1)(X_1Z_2 + X_2Z_1) - 3B(X_1Z_2-X_2Z_1)Z_1Z_2, \\[0.2cm]
    S_2    = \ & -3X_1X_2(X_1Y_2-X_2Y_1) - Y_1Y_2(Y_1Z_2-Y_2Z_1) - A(X_1Y_2-X_2Y_1)Z_1Z_2 \\
            & + A(Y_1Z_2-Y_2Z_1)(X_1Z_2 + X_2Z_1) + 3B(Y_1Z_2-Y_2Z_1)Z_1Z_2, \\[0.2cm]
    S_3    = \ & 3X_1X_2(X_1Z_2-X_2Z_1) - (Y_1Z_2-Y_2Z_1)(Y_1Z_2 + Y_2Z_1) \\
            & + A(X_1Z_2-X_2Z_1)Z_1Z_2,
\end{align*}
and
\begin{align*}
    T_1    = \ & Y_1Y_2(X_1Y_2 + X_2Y_1) - AX_1X_2(Y_1Z_2 + Y_2Z_1) \\
                        & - A(X_1Y_2 + X_2Y_1)(X_1Z_2 + X_2Z_1) - 3B(X_1Y_2 + X_2Y_1)Z_1Z_2 \hspace{1.6cm} \\
                        & - 3B(X_1Z_2 + X_2Z_1)(Y_1Z_2 + Y_2Z_1) + A^2(Y_1Z_2 + Y_2Z_1)Z_1Z_2, \\[0.2cm]
    T_2    = \ & Y_1^2Y_2^2 + 3AX_1^2X_2^2 + 9BX_1X_2(X_1Z_2 + X_2Z_1) \\
                        & - A^2X_1Z_2(X_1Z_2 + 2X_2Z_1) - A^2X_2Z_1(2X_1Z_2 + X_2Z_1) \\
                        & - 3ABZ_1Z_2(X_1Z_2 + X_2Z_1) - (A^3 + 9B^2)Z_1^2Z_2^2, \\[0.2cm]
    T_3    = \ & 3X_1X_2(X_1Y_2 + X_2Y_1) + Y_1Y_2(Y_1Z_2 + Y_2Z_1) \\
                        & + A(X_1Y_2 + X_2Y_1)Z_1Z_2 + A(X_1Z_2 + X_2Z_1)(Y_1Z_2 + Y_2Z_1) \\
                        & + 3B(Y_1Z_2 + Y_2Z_1)Z_1Z_2.
\end{align*}
A compact and efficient way of computing the latter addition law may be found in \cite[Lemma 2.1]{ELoops}.
Similar concise formulas over any characteristics were established in \cite[Proposition 3.2]{InvTau}.


\section{Elliptic curves over $\Z/N\Z$}
\label{S:3}

Let $N \in \Z_{\geq 2}$ be an integer. Hereafter we consider elliptic curves defined over the ring $R = \Z/N\Z$, which satisfies the Condition \ref{R-Condition}.
More generally, in \cite{ECNTA} this condition has been proved to hold for every ring with a finite number of maximal ideals.
Here we show that $\Z/N\Z$ underlies a condition that is even stronger than Condition \ref{R-Condition}.

\begin{lemma} \label{lem:primitiveoverN}
Let $N \in \Z_{\geq 2}$ be an integer and $A$ be a matrix over $\Z/N\Z$ whose entries are primitive, then there exists a linear combination of the columns of $A$ that is primitive.
In particular, $R = \Z/N\Z$ satisfies Condition \ref{R-Condition}.
\end{lemma}
\proof Let $A = (c_1 | c_2 | \dots | c_m)$ be the columns of the considered matrix.
Since $A$ is primitive, for every prime $p | N$ there are coefficients $\a^{(p)}_1, \dots, \a^{(p)}_m \in \Z/p\Z$ such that the vector
\[
    v^{(p)} = \sum_{i = 1}^m \a^{(p)}_i c_i
\]
is primitive over $\Z/p\Z$.
By the Chinese Reminder Theorem we may find integers $\b_1, \dots, \b_m \in \Z$ solving, for every prime divisor $p$ of $N$, the congruence system
\[
    \b_i \equiv \a^{(p)}_i \bmod p.
\]
Therefore, $\sum_{i = 1}^m \b_i c_i$ is easily seen to be a primitive combination of the columns of $A$.
\endproof

We now recall how the group of points of an elliptic curve over $\Z/N\Z$ can be described by the curve projections over the $p$-components of this ring, with $p$ ranging among the prime divisors of $N$.

\begin{proposition}[\cite{Washington}, Corollary 2.32] \label{prop:WashingtonSplit}
Let $N_1, N_2$ be coprime integers and let $A,B \in \Z$ such that $\Delta_{A,B} \in (\Z/N_1N_2\Z)^*$. Then the canonical projections induce a group isomorphism
\[
     E_{A,B}(\Z/N_1N_2\Z) \simeq E_{A,B}(\Z/N_1\Z) \oplus E_{A,B}(\Z/N_2\Z).
\]
\end{proposition}

\noindent Thus, it is sufficient to study the structure of elliptic curves $E_{A,B}(\Z/p^e\Z)$ for any prime $p$ and positive integer $e$, which is the main goal of this section.
We begin by noticing that the points $P \in E = E_{A,B}(\Z/p^e\Z)$ of such curves have prescribed representatives:

\begin{itemize}
    \item If $P \in E^a$, then there are $X,Y \in \Z/p^e\Z$ such that
    
    \[
        P = (X:Y:1).
    \]
    \item If $P \in E^{\infty}$, then there are $X,Z \in p(\Z/p^e\Z)$ such that
    
    \[
        P = (X:1:Z).
    \]    
\end{itemize}


The size of these curves is known, as reported in the next lemma.

\begin{lemma}[\cite{ECNTA}, Section 4] \label{lem:cardinality}
Let $p$ be a prime, $e \in \Z_{\geq 1}$ and 
\[
    \pi : E_{A,B}(\Z/p^e\Z) \to E_{A,B}(\F_p)
\]
be the canonical projection.
Then for every $P \in E_{A,B}(\F_p)$ we have
\[
    |\pi^{-1}(P)| = p^{e-1}.
\]
In particular:
\begin{itemize}
    \item the size of the curve is $|E_{A,B}(\Z/p^e\Z)| = p^{e-1} |E_{A,B}(\F_p)|$,
    \item $\ker \pi$ is a subgroup of $E_{A,B}(\Z/p^e\Z)$, whose size is $p^{e-1}$.
\end{itemize}
\end{lemma}

The coordinates of points at infinity
satisfy the following relation, which we prove by adapting the idea of expansion around $\0$ \cite[Chapter IV]{Silverman}.

\begin{proposition} \label{prop:InfPoints}
Let p be a prime, $e \in \Z_{\geq 1}$ and $E = E_{A,B}(\Z/p^e\Z)$.
There is a polynomial $\ff \in \Z[x]$ of degree at most $e-1$ such that for every $P \in E^{\infty}$ there is $X \in p(\Z/p^e\Z)$ satisfying
\[
    P = \big(X:1:\ff(X)\big).
\]
Moreover, we have
\[
    \ff(X) \equiv X^3 + AX^7 + BX^9 \bmod{p^{10}}.
\]
\end{proposition}

\proof Since $P \in E^{\infty}$, it may be represented in the form $(X:1:Z)$, with $X,Z \in p(\Z/p^e\Z)$ that satisfy
\[
    Z \equiv X^3 + AXZ^2 + BZ^3 \bmod{p^e}.
\]
We recursively define the following sequence of polynomials in $\Z[x,z]$:
\[
    F_0(x,z) = x^3 + Axz^2 + Bz^3, \qquad \forall \ i \in \Z_{\geq 1} : \ F_i(x,z) = F_{i-1}\big(x,F_0(x,z)\big).
\]
It is easy to see by induction on $i \in \Z_{\geq 0}$ that this sequence satisfies 
\[
    Z \equiv F_i(X,Z) \bmod{p^e}.
\]
Moreover, every $F_i$ for $i \in \Z_{\geq 1}$ is obtained from $F_{i-1}$ by substituting all the occurrences of $z$ with $F_0(x,z)$, which contains only terms of degree $3$, hence the total degree of terms involving $z$ in $F_i$ is strictly increasing while increasing $i$.
This means that there exist $M \in \Z_{\geq 0}$ and $g \in \Z[x,z]$ such that
\[
    F_M(x,z) = \ff(x) + g(x,z), \qquad \text{ with } \begin{cases}\deg(g) \geq e,\\ \deg(\ff) < e.\end{cases}
\]
Since both $X$ and $Z$ are divisible by $p$, this implies
\[
    Z \equiv F_M(X,Z) \equiv \ff(X) \bmod{p^e},
\]
so that $\ff \in \Z[x]$ is the required polynomial.
A direct computation shows that
\[
F_3 = x^3 + Ax^7 + Bx^9 + (\text{terms of degree} \geq 11),
\]
which proves the moreover part.
\endproof

\begin{remark}
Although finite local rings are complete with respect to the topology induced by their maximal ideal, they may well not be domains (e.g. $\Z/N\Z$).
For this reason, we found it appropriate to explicitly compute $\ff$ instead of considering the truncation to the correct exponent of the classical series \cite[Chapter IV]{Silverman}.
\end{remark}

To simplify the exposition, for any $X \in \Z/p^e\Z$ and any positive integer $t$ we write $p^t | X$ or $X \equiv 0 \bmod p^t$ in place of the more precise $X \in p^t(\Z/p^e\Z)$.
In the same spirit, we assign a $p$-adic valuation to any $X \in \Z/p^e\Z$ by writing
\[
    \v_p(X) = \begin{cases} t & \text{if } X \in p^t(\Z/p^e\Z) \setminus p^{t+1}(\Z/p^e\Z),\\
    e & \text{if } X = 0.
    \end{cases}
\]

From Proposition \ref{prop:InfPoints} it is possible to derive a description of the first-order approximation of the sum of two points at infinity.

\begin{proposition} \label{prop:InfSum}
Let $p$ be a prime, $e \in \Z_{\geq 1}$, $E = E_{A,B}(\Z/p^e\Z)$ and $\ff \in \Z[x]$ be the polynomial arising from $E$ as in Proposition \ref{prop:InfPoints}. Let also
\[
    P_1 = \big(X_1:1:\ff(X_1)\big),\ P_2 = \big(X_2:1:\ff(X_2)\big) \in E^{\infty}
\]
with $e_1 = \v_p(X_1)$ and $e_2 = \v_p(X_2)$.
Then
\[
    P_1+P_2 = \big(X_3:1:\ff(X_3)\big),
    \quad \text{where} \quad X_3 \equiv X_1 + X_2 \bmod{p^{5\min\{e_1,e_2\}}}.
\]
\end{proposition}
\proof As $\pi$ is a group homomorphism, $P_1+P_2$ lies itself in $E^{\infty}$, which implies [\cite{BosmaLenstra}, Theorem 2] that these points are never exceptional for the addition law $+_{(0:1:0)}$ corresponding to $(0:1:0)$.
A straightforward computation with $+_{(0:1:0)}$ shows that, modulo monomials in $X_1$ and $X_2$ of total degree at least $5$ (i.e. modulo $p^{5\min\{e_1,e_2\}}$), we have
\[
    P_1 + P_2 = \big(X_1 + X_2 : 1 + 3AX_1^2X_2^2 : (X_1 + X_2)^3\big),
\]
which is equal to $\big(X_1 + X_2 : 1 : (X_1 + X_2)^3\big)$ as we verify by multiplying its entries by $1-3AX_1^2X_2^2 \in (\Z/p^{5\min\{e_1,e_2\}}\Z)^*$.
\endproof


We can now prove that the infinity group is cyclic, which provides a structure theorem for elliptic curves over $\Z/N\Z$.

\begin{theorem} \label{thm:structure}
Let $p$ be a prime, $e \in \Z_{\geq 1}$ and $\ff \in \Z[x]$ be the polynomial arising from $E_{A,B}(\Z/p^e\Z)$ as in Proposition \ref{prop:InfPoints}. Then
\[
    0 \to \big\langle \big(p:1:\ff(p)\big) \big\rangle
    \stackrel{\iota}{\hookrightarrow}
     E_{A,B}(\Z/p^e\Z) \stackrel{\pi}{\to} E_{A,B}(\F_p) \to 0.
\]
is a short exact sequence of groups.
\end{theorem}
\proof We know that the canonical projection $\pi : E_{A,B}(\Z/p^{e}\Z) \twoheadrightarrow E_{A,B}(\F_p)$ is a surjective group homomorphism and that $|\ker \pi| = p^{e-1}$ by Lemma \ref{lem:cardinality}.
Thus, it is sufficient to prove that $P = \big(p:1:\ff(p)\big) \in \ker \pi$ has order $p^{e-1}$.
Since $P$ lies over $\0 \in E_{A,B}(\F_p)$, then its order is a power of $p$ ($\ker \pi$ is a $p$-group).
We prove by induction on $0 \leq \e \leq e-1$ that
\[
    p^{\e} P = \big(X:1:\ff(X)\big) \quad \text{ and } \quad \v_p(X) = \e+1.
\]
In particular, the minimal $\e$ such that $X \equiv 0 \bmod p^e$ is $\e = e-1$.\\
$[\e = 0]$ It is trivially seen that
\[
    p^0 P = \big(p:1:\ff(p)\big) \quad \text{ and } \quad \v_p(p) = 1.
\]
$[\e \rightarrow \e+1]$ By the inductive hypothesis we know that
\[
    p^{\e+1} P = p (p^\e P) = p \big(X:1:\ff(X)\big) \quad \text{ and } \quad \v_p(X) = \e+1.
\]
By Proposition \ref{prop:InfSum} and induction on $\a \in \{1, \dots, p-1\}$ we have
\[
    \big(X:1:\ff(X)\big) + \big(\a X:1:\ff(\a X)\big) = \big(X_2:1:\ff(X_2)\big),
\]
with
\[
    X_2 \equiv (\a+1)X \bmod p^{5(\e+1)}.
\]
Thus, by specializing the above result for $\a = p-1$, the $p$-adic valuation of the first component of $p\big(X:1:\ff(X)\big)$ is proved to be $\v_p(X)+1 = \e+2$.
\endproof

The above theorem shows that the infinity part of any elliptic curve over $\Z/p^e\Z$ is a $\Z/p^e\Z$-torsor with respect to the standard multiplication action.
This agrees with \cite[Proposition 10.3]{ELoops} and it is sufficient to determine the group structure of these curves when their projection is not anomalous.

\begin{corollary} \label{cor:structure}
Let $p$ be a prime, $e \in \Z_{\geq 1}$ and $E_{A,B}(\Z/p^e\Z)$ be an elliptic curve such that $|E_{A,B}(\F_p)| \neq p$. Then
\[
    E_{A,B}(\Z/p^e\Z) \simeq E_{A,B}(\F_p) \oplus \Z/p^{e-1}\Z.
\]
\end{corollary}
\proof It is sufficient to show that the short exact sequence of Theorem \ref{thm:structure} splits, which by the Splitting Lemma amounts to proving that it is left split.
Since $q = |E_{A,B}(\F_p)| \neq p$ is in the Hasse bound of $p$, then $(p,q) = 1$, which implies that there is $k \in \Z$ satisfying
\[
    \begin{cases}
    k \equiv 1 \bmod p^{e-1},\\
    k \equiv 0 \bmod q.
    \end{cases}
\]
By Theorem \ref{thm:structure} we have $E_{A,B}^{\infty}(\Z/p^e\Z) = \pi^{-1}(\0) = \big\langle \big(p:1:\ff(p)\big) \big\rangle$. Thus, since $k \equiv 0 \bmod q$, the map
\[
    E_{A,B}(\Z/p^{e}\Z) \stackrel{\cdot k}{\to} \big\langle \big(p:1:\ff(p)\big) \big\rangle
\]
is a well-defined group homomorphism.
Moreover, since $k \equiv 1 \bmod p^{e-1}$, the cyclic group $\big\langle \big(p:1:\ff(p)\big) \big\rangle$ is fixed under this map, hence the multiplication-by-$k$ is a left section for the considered sequence.
\endproof

When $e$ is small, an explicit group isomorphism may also be exhibited.

\begin{proposition} \label{Prop:ExplicitIso}
Let $p$ be a prime, $1 \leq e \leq 5$ be an integer, $E_{A,B}(\Z/p^e\Z)$ be an elliptic curve and $q = |E_{A,B}(\F_p)|$ be the size of its projected curve. Then
\begin{align*}
    \Phi : E_{A,B}(\Z/p^e\Z) &\to E_{A,B}(\F_p) \oplus \Z/p^{e-1}\Z,\\
    P &\mapsto \left( \pi(P), \frac{1}{p} \frac{(qP)_x}{(qP)_y} \right),
\end{align*}
is a well-defined group homomorphism.
Moreover, if $q \neq p$ then $\Phi$ is a group isomorphism.
\end{proposition}
\proof It is easy to see that $\Phi(P)$ does not depend on the projective representative of $P$.
Moreover, as $\pi$ is a group homomorphism we have
\[
    \pi(qP) = q \pi(P) = \0 \in E_{A,B}(\F_p).
\]
Hence, by Proposition \ref{prop:InfPoints} we have $qP = \big(X:1:\ff(X)\big)$ with $X \in p(\Z/p^e\Z)$.
Therefore, $\frac{(qP)_x}{(qP)_y} \in p(\Z/p^{e}\Z)$, which is canonically isomorphic to $\Z/p^{e-1}\Z$.
Thus, $\Phi$ is a well-defined map between groups having, by Lemma \ref{lem:cardinality}, the same size.
It also respects the sum, as for every pair $P_1, P_2 \in E_{A,B}(\Z/p^{e}\Z)$ we compute
\[
    \Phi(P_1) + \Phi(P_2) = \left( \pi(P_1 + P_2), \frac{1}{p} \left(\frac{(qP_1)_x}{(qP_1)_y} +  \frac{(qP_2)_x}{(qP_2)_y}\right) \right),
\]
and since $e \leq 5\min\{\v_p\big((qP_1)_x\big), \v_p\big((qP_2)_x\big)\}$, then by Proposition \ref{prop:InfSum} we have
\[
    \frac{(qP_1)_x}{(qP_1)_y} + \frac{(qP_2)_x}{(qP_2)_y} = \frac{(qP_1 + qP_2)_x}{(qP_1 + qP_2)_y} = \frac{\big(q(P_1 + P_2)\big)_x}{\big(q(P_1 + P_2)\big)_y}.
\]

As for the moreover part, it is sufficient to prove that $\ker \Phi = \{\0\}$ when $q \neq p$.
Let $\Phi(P) = (\0, 0)$, then there exists $X \in p(\Z/p^e\Z)$ such that $P = \big(X:1:\ff(X)\big)$ and
\[
    \frac{qX}{p} \equiv \frac{(qP)_x}{p} \equiv 0 \bmod{p^{e-1}}.
\]
Since $q$ lies in the Hasse interval of $p$, then $q \neq p$ implies $(p,q) = 1$ and we conclude that $X \equiv 0 \bmod{p^e}$, hence the kernel of $\Phi$ is trivial.
\endproof

When the restricted curve $E_{A,B}(\F_p)$ is anomalous two different scenarios may occur.
By Theorem \ref{thm:structure} the curve $E_{A,B}(\Z/p^e\Z)$ is guaranteed to contain a cyclic subgroup of order $p^{e-1}$, therefore it may be either cyclic
\begin{equation} \tag{Cyclic}
    E_{A,B}(\Z/p^e\Z) \simeq \Z/p^e\Z,
\end{equation}
or split, i.e.
\begin{equation} \tag{Split}
    E_{A,B}(\Z/p^e\Z) \simeq \F_p \oplus \Z/p^{e-1}\Z.
\end{equation}

Even if the cyclic case occurs over $\Z/p^e\Z$ with an overwhelming probability ($\frac{p-1}{p}$), both may take place. As an instance, one may check that

\begin{small}
\[
    E_{7,3}(\Z/13^2\Z) \simeq \langle (0:61:1) \rangle, \quad \text{while} \quad E_{1,6}(\Z/13^2\Z) \simeq \langle (2:4:1) \rangle \oplus \langle (13:1:0) \rangle.
\]
\end{small}
The above discussion leads to the classification theorem.

\begin{theorem} \label{thm:classification}
Let $N$ be a positive integer and let $A,B$ be integers such that $\Delta_{A,B}$ is coprime to $N$. Then we have
\[
    E_{A,B}(\Z/N\Z) \simeq \bigoplus_{\substack{p | N\\ |E_{A,B}(\F_p)| \neq p}} E_{A,B}(\F_p) \oplus \Z/p^{\v_p(N)-1}\Z \oplus \bigoplus_{\substack{p | N\\ |E_{A,B}(\F_p)| = p}} G_p,
\]
where every $G_p$ may be either $\Z/p^{\v_p(N)}\Z$ or $\F_p \oplus \Z/p^{\v_p(N)-1}\Z$.
\end{theorem}
\proof By Proposition \ref{prop:WashingtonSplit} we know that
\[
    E_{A,B}(\Z/N\Z) \simeq \bigoplus_{p | N} E_{A,B}(\Z/p^{\v_p(N)}\Z).
\]
By Corollary \ref{cor:structure}, for every $p$ such that $E_{A,B}(\F_p)$ is not anomalous we have
\[
    E_{A,B}(\Z/p^e\Z) \simeq E_{A,B}(\F_p) \oplus \Z/p^{\v_p(N)-1}\Z.
\]
On the other side, we have seen that
\[
    G_p = \F_p \oplus \Z/p^{\v_p(N)-1}\Z \qquad \text{or} \qquad G_p = \Z/p^{\v_p(N)}\Z
\]
may both occur as group structure of $E_{A,B}(\Z/p^{\v_p(N)}\Z)$ when $E_{A,B}(\F_p)$ is anomalous, which completes the study cases.
\endproof

\begin{remark}
Given a finite collection of elliptic curves $\{E_{A_l,B_l}(R_l)\}_{1 \leq l \leq k}$, we may define an elliptic curve over their product ring $\prod_{l=1}^k R_l$ with the componentwise operation and \cite[Section 4]{ECNTA} we have
\[
    E_{(A_1,\dots,A_k),(B_1,\dots,B_k)}\left(\prod_{l=1}^k R_l\right) \simeq \prod_{l=1}^k E_{A_l,B_l}(R_l).
\]
Thus, Theorem \ref{thm:classification} provides the group structures of every elliptic curve defined over a ring isomorphic to a finite product of integer residue rings.
\end{remark}

\begin{remark}
We notice that Theorem \ref{thm:structure} really relies on the behavior of elliptic curves over $\Z/p^e\Z$.
Let us consider another local ring, namely $R = \F_5[x]/(x^4)$, and let $\e$ be a generator of its maximal ideal.
It defines a canonical projection
\[
    R \to \F_5, \quad X_0 + X_1\e + X_2\e^2 + X_3\e^3 \mapsto X_0,
\]
so we have an elliptic curve $E_{1,2}(R)$ defined as in Section \ref{subsec:ECoverR}, together with a canonical projection onto $E_{1,2}(\F_5)$.

This curve may appear similar to $E_{1,2}(\Z/5^4\Z)$ at first glance, but one can directly verify that $E_{1,2}(R)$ is given by
\[
     \langle (2\e^3 + \e: 1: \e^3) \rangle \oplus \langle ( 3\e^3 + 3\e^2 + 2\e: 1: 3\e^3) \rangle \oplus \langle (\e^3 + \e + 3:\e^3 + 3\e^2 + 4\e + 3:1) \rangle,
\]
so that $E_{1,2}(R) \simeq \Z/5\Z \oplus \Z/5\Z \oplus \Z/35\Z$. This is due to the different structure of the infinity parts, as $E^{\infty}_{1,2}(R) \simeq (\Z/5\Z)^{\oplus 3}$, while $E^{\infty}_{1,2}(\Z/5^4\Z) \simeq \Z/5^3\Z$ as prescribed by our previous results.
A detailed study of the latter type of rings may be found in \cite{InvTau}.
\end{remark}


\section{Rank of $p$-groups from elliptic curves} \label{S:4}

We know that groups arising from elliptic curves defined over finite fields have prescribed constraints \cite{Ruck,Voloch}, e.g. their rank cannot exceed $2$. 
This restriction can be relaxed for curves defined over $\Z/N\Z$, as their rank may be arbitrarily large, but it may still be bounded in terms of the number of primes inside a Hasse interval.

\begin{definition}[$\H_p$]
Given an integer $p \in \Z$, we define
\[
    \H_p = |\{q \in \Z \ | \ q \text{ is prime and } p+1-2\sqrt{p} \leq q \leq p+1+2\sqrt{p}\}|.
\]
\end{definition}

The following result provides a sharp bound on the rank that elliptic curves over $\Z/N\Z$ may have if their point group are $p$-groups, which in particular shows that there are infinitely many groups that cannot arise as a point group for an elliptic curve over an integer residue ring.

\begin{proposition} \label{prop:SpecialRk}
Let $p \geq 5$ be a prime, $N \in \Z_{\geq 2}$ and $E=E_{A,B}(\Z/N\Z)$ be an elliptic curve that is a $p$-group.
Then, by defining
\[
    \chi_p = \begin{cases} 2 &\text{if there is a prime $q$ such that } E_{A,B}(\F_q) \simeq \F_p \oplus \F_p,\\
    0 &\text{otherwise},
    \end{cases}
\]
we have
\[
    \rk(E) \leq \H_p + \chi_p + 1.
\]
\end{proposition}
\proof By Theorem \ref{thm:classification} we have
\[
    E \simeq \bigoplus_{\substack{q | N\\ |E_{A,B}(\F_q)| \neq q}} E_{A,B}(\F_q) \oplus \Z/q^{\v_q(N)-1}\Z \oplus \bigoplus_{\substack{q | N\\ |E_{A,B}(\F_q)| = q}} G_q,
\]
where every $G_q$ may be either $\Z/q^{\v_q(N)}\Z$ or $\F_q \oplus \Z/q^{\v_q(N)-1}\Z$.
It is easy to see that $G_q$ is a $p$-group only if $q = p$, hence we have
\[
    \rk \bigoplus_{\substack{q | N\\ |E_{A,B}(\F_q)| = q}} G_q \leq 2.
\]
Similarly, we notice that $\Z/q^{\v_q(N)-1}\Z$ is a $p$-group only if $q = p$, but $E_{A,B}(\F_p)$ is a $p$-group if and only if $|E_{A,B}(\F_p)| = p$.
Thus, we have
\[
    \bigoplus_{\substack{q | N\\ |E_{A,B}(\F_q)| \neq q}} E_{A,B}(\F_q) \oplus \Z/q^{\v_q(N)-1}\Z \simeq \bigoplus_{\substack{q | N\\ q \neq p}} E_{A,B}(\F_q).
\]
Moreover, since the rank of $E_{A,B}(\F_q)$ is at most $2$ \cite[Theorem 4.1]{Washington}, then it is a $p$-group only if
\[
    \text{either} \qquad E_{A,B}(\F_q) \simeq \F_{p} \qquad \text{or} \qquad E_{A,B}(\F_q) \simeq \F_{p} \oplus \F_{p}.
\]
Since the Hasse buond over a prime field is full, then $E_{A,B}(\F_q)$ may be isomorphic to $\F_p$ for every prime $q$ inside the Hasse interval of $p$.

On the other side, by \cite[Prop.4.16]{Washington} we know that $E_{A,B}(\F_q) \simeq \F_{p} \oplus \F_{p}$ may occur only if 
\[
    q \in \{p^2+1, p^2\pm p+1, p^2\pm 2p+1\}.
\]
However, both $p$ and $q$ are odd primes, hence only $q = p^2\pm p+1$ may occur.
Furthermore, since $p > 3$, it is easy to see that either $3 | p^2+p+1$ or $3 | p^2-p+1$, therefore only one of them can be prime.
We conclude that there is at most one prime $q$ such that $E_{A,B}(\F_q) \simeq \F_{p} \oplus \F_{p}$, so that
\[
    \rk \bigoplus_{\substack{q | N\\ q \neq p}} E_{A,B}(\F_q) \leq (\H_p - 1) + \chi_p.
\]
Collecting the above rank bounds, the statement follows.
\endproof

\begin{example}
Let $p=11$.
None of $11^2\pm 11+1$ is prime, then we have $\chi_{11} = 0$, therefore by Proposition \ref{prop:SpecialRk}, regardless of $N \in \Z_{\geq 2}$, the rank of any elliptic curve over $\Z/N\Z$ that is a $11$-group is bounded by $\H_{11}+1 = 5$.
We also notice that this bound is sharp, as
\[
    E_{167707,21664}(\Z/ 187187\Z) \simeq \F_{11} \oplus \F_{11} \oplus \F_{11} \oplus \F_{11} \oplus \F_{11}.
\]
\end{example}

\begin{example}
Let $p=13$.
We notice that $13^2-13+1=157$ is prime and
\[
    E_{0,15}(\F_{157}) \simeq \F_{13} \oplus \F_{13},
\]
therefore we have $\chi_{13} = 2$.
By means of Proposition \ref{prop:SpecialRk} we know that any elliptic curve over $\Z/N\Z$ that is a $13$-group has rank bounded by $\H_{13}+3 = 8$.
We notice again that this bound is sharp, as
\[
    E_{63707931,239467091}(\Z/ 659902243\Z) \simeq (\F_{13})^{\oplus 8}.
\]
\end{example}


\section{Another isomorphism attack to anomalous ECDLP} \label{S:5}

Given an additive group $G$ and a base element $g \in G$, the discrete logarithm problem (DLP) on $G$ consists of computing for any given $h \in G$ an integer $N$, if existent, such that $h = N \cdot g = g + g + \dots + g$.
When $G$ is the point group of an elliptic curve (ECDLP), this problem is known to be computationally feasible only in special cases, such as the anomalous ones \cite{Semaev,SatohAraki,Smart}.

From the knowledge of the group structure provided by Theorem \ref{thm:structure}, we have another way for efficiently solving the ECDLP on anomalous curves by employing any cyclic curve that projects onto it.

\begin{proposition} \label{prop:cyclicStruct}
Let $p$ be a prime, $e\in \Z_{\geq 2}$ and $E = E_{A,B}(\Z/p^{e}\Z)$ be an elliptic curve, whose point group is cyclic of order $p^e$. Then the map
\begin{align*}
    \Theta : E &\to \F_p,\\
    P &\mapsto \frac{1}{p^{e-1}} \frac{(p^{e-1}P)_x}{(p^{e-1}P)_y},
\end{align*}
is a well-defined surjective group homomorphism, whose kernel is
\[
    \ker \Theta = \big\langle \big(p:1:\ff(p)\big) \big\rangle.
\]
\end{proposition}

\proof For every $P \in E$ the point $p^{e-1}P$ is a $p$-torsion point of $E$, hence
\[
    p^{e-1}P = \big(X:1:\ff(X)\big), \qquad \text{ with }\ \v_p(X) \geq e-1,
\]
therefore $\Theta(P) = \frac{X}{p^{e-1}} \in \F_p$ is well-defined.
Let $G \in E$ be a generator of the point group of $E$, then for every integer $m \in \Z$ we have
\[
    p^{e-1}m G = m \big(X : 1 : \ff(X)\big) = \big(m X : 1 : \ff(mX)\big),
\]
where the last equality follows from Proposition \ref{prop:InfSum}, as for every $e \geq 2$ the point $p^{e-1}G$ lies in $\langle (p^{e-1}:1:0 ) \rangle$. Thus, $\Theta( m G ) = m \Theta( G )$, so that $\Theta$ is a group homomorphism.
Moreover, from the above equation it follows that
\[
    \ker \Theta = \{ mp\ G \ | \ m \in \Z \} = \big\langle \big(p:1:\ff(p)\big) \big\rangle.
\]
By comparing the size of these groups, the surjectivity follows.
\endproof

From the above proposition, the discrete logarithm over anomalous curves may be immediately recovered.

\begin{corollary} \label{cor:ANDLOG}
Let $p$ be a prime, $e\in \Z_{\geq 2}$ and $E_{A,B}(\Z/p^{e}\Z)$ be an elliptic curve, whose point group is cyclic of order $p^e$. Then the map
\[
    \Theta \circ \pi^{-1} : E_{A,B}(\F_p) \to \F_p
\]
is a well-defined group isomorphism.
\end{corollary}
\proof By Theorem \ref{thm:structure} the canonical projection induces a group isomorphism $E_{A,B}(\Z/p^{e}\Z)/\langle (p:1:\ff(p)) \rangle \simeq E(\F_p)$, whereas the map $\Theta$ arisen from Proposition \ref{prop:cyclicStruct} induces a group isomorphism $E_{A,B}(\Z/p^{e}\Z)/\langle (p:1:\ff(p)) \rangle \simeq \F_p$.
By composing those isomorphisms, the result follows.
\endproof

Finding any lift of a given point is computationally costless, therefore the complexity of the isomorphism attack given by Corollary \ref{cor:ANDLOG} only depends on the cost of computing $\Theta$, which is $O(\log p)$.
This approach is not faster than previously known attacks to the same family of curves, but it has the advantage of involving only finite precision objects.

\begin{example}
Let us consider an anomalous curve as constructed in \cite{HugeAnomalous}:
\begin{align*}
    p &= 730750818665451459112596905638433048232067471723,\\
    A &= 425706413842211054102700238164133538302169176474,\\
    B &= 203362936548826936673264444982866339953265530166.
\end{align*}
We consider on $E_{A,B}(\F_p)$ the points 
\begin{align*}
    P &= ( 1 : 310536468939899693718962354338996655381367569020 : 1 ), \\
    Q &= ( 3 : 38292783053156441019740319553956376819943854515 : 1 ).
\end{align*}
To find their discrete logarithm it is sufficient to compute any lifts, such as
\[
    P^{\uparrow} = ( 1 : P_y + \a p : 1 ),\ Q^{\uparrow} = ( 3 : Q_y + \b p : 1 ) \in E_{A,B}(\Z/p^2\Z),
\]
where
\[
    \a = \frac{1+A+B-P_y^2}{2pP_y} \bmod p^2,\ \b = \frac{27+3A+B-Q_y^2}{2pQ_y} \bmod p^2,
\]
and to apply them the group homomorphism $\Theta$ of Proposition \ref{prop:cyclicStruct}:
\begin{align*}
    \Theta(P^{\uparrow}) &= 343088892565802863386490109374548044078624360215, \\
    \Theta(Q^{\uparrow}) &= 470974712001084540433398653921983741661987449793.
\end{align*}
This way we get the discrete logarithm $N$ such that $Q = N \cdot P$ as

\begin{small}\[
    N = \frac{\Theta(Q^{\uparrow})}{\Theta(P^{\uparrow})} \bmod p = 113690975836469390483838646646828917131453128585.
\]
\end{small}

We remark that such a discrete logarithm would be unfeasible to be computed with generic logarithm techniques, as one can directly verify that the \texttt{Log} routine of Magma \cite{Magma} does not terminate in a reasonable time. 
\end{example}


\section{Conclusions and open problems} \label{S:6}

In this work, we have provided the classification of groups arising from elliptic curves over $\Z/N\Z$ and exploited it for obtaining a bound for their rank and an attack on the ECDLP over anomalous elliptic curves.

The key ingredient is Theorem \ref{thm:structure}, which might still hold for more general classes of rings, even though the kernel generator may be less explicit.
Finding other instances or even classifying all the rings over which the infinity group is cyclic is still an open line of research.

From a cryptographic perspective, Theorem \ref{thm:structure} shows that the difficulty of the ECDLP depends on the difficulty of the same problem over the base field and in the group of points at infinity. Whenever these two groups are linked (as in the case of the anomalous curves), the discrete logarithm on one group may be read from the other.

Finally, in this work we only considered genus-$1$ curves for their theoretical and historical relevance, but it is reasonable to ask which other abelian varieties admit such an extension to $\Z/N\Z$ and, when it is the case, if analogous group decompositions over these rings hold.
%
%
%
 \bibliographystyle{splncs04}

\begin{thebibliography}{8}
\bibitem{BPS}
W. D. Banks, F. Pappalardi, and I. E. Shparlinski,
\emph{On Group Structures Realized by Elliptic Curves over Arbitrary Finite Fields},
Experiment. Math. 21, 2012, pp. 11--25.

\bibitem{BosmaFact}
W. Bosma,
\emph{Primality testing using elliptic curves},
Math. Instituut, Univ. Amsterdam, volume Tech. Rep. 85--12, 1985.

\bibitem{BosmaLenstra}
W. Bosma, and H. W. Lenstra,
\emph{Complete Systems of Two Addition Laws for Elliptic Curves},
J. Number Theory 53, 1995, pp. 229--240.

\bibitem{Magma}
W. Bosma, J. Cannon, and Catherine Playoust,
\emph{The Magma algebra system. I. The user language},
J. Symbolic Comput. 24, 1997, pp. 235--265.

\bibitem{Modularity}
C. Breuil, B. Conrad, F. Diamond, and R. Taylor,
\emph{On the modularity of elliptic curves over Q: wild 3-adic exercises},
J. Amer. Math. Soc. 14, 2001, pp. 843--939.

\bibitem{Brown}
W. C. Brown,
\emph{Matrices over commutative rings},
Marcel Dekker, 1986.

\bibitem{Husemoller}
D. Husem\"oller,
\emph{Elliptic Curves},
Grad. Texts in Math. 111. Springer-Verlag, 1987.

\bibitem{InvTau}
R. Invernizzi, and D. Taufer,
\emph{Multiplication polynomials for elliptic curves over finite local rings},
ACM’s International Conference Proceedings Series (ISSAC 2023), 2023, pp. 335--344. 

\bibitem{ECDSA}
D. Johnson, A. Menezes, and S. Vanstone,
\emph{The Elliptic Curve Digital Signature Algorithm (ECDSA)},
Int. J. Inf. Secur. 1, 2001, pp. 36--63.

\bibitem{Koblitz}
N. Koblitz,
\emph{Elliptic curve cryptosystems},
Math. Comp. 48, 1987, pp. 203--209.

\bibitem{GrpGenerators}
D. R. Kohel, and I. E. Shparlinski,
\emph{On Exponential Sums and Group Generators for Elliptic Curves over Finite Fields},
Lecture Notes in Comput. Sci. 21, 2000, pp. 395--404.

\bibitem{CryptoZNZ}
K. Koyama, U. M. Maurer, T. Okamoto, and S. A. Vanstone,
\emph{New Public-Key Schemes Based on Elliptic Curves over the Ring Z$_n$},
Adv. Cryptology 576, 1991, pp. 252--266.

\bibitem{LanRup85}
H. Lange, and W. Ruppert,
Complete systems of addition laws on abelian varieties,
Invent. Math. 79, 1985, pp. 603--610.

\bibitem{ECNTA}
H. W. Lenstra,
\emph{Elliptic curves and number-theoretic algorithms},
Proc. of the International Congress of Mathematicians, 1986, pp. 99--120.

\bibitem{LenstraFact}
H. W. Lenstra,
\emph{Factoring integers with elliptic curves},
Ann. of Math. 126, 1987, pp. 649--673.

\bibitem{HugeAnomalous}
F. Leprévost, J. Monnerat, S. Varrette, and S. Vaudenay,
\emph{Generating Anomalous Elliptic Curves},
Inform. Process. Lett. 93, 2005, pp. 225--230. 

\bibitem{Merel}
L. Merel,
\emph{Bornes pour la torsion des courbes elliptiques sur les corps de
nombres},
Invent. Math. 124, 1996, pp. 437--449.

\bibitem{Miller}
V. S. Miller,
\emph{Use of elliptic curves in cryptography},
Adv. Cryptol. 218, 1985, pp. 417--426.

\bibitem{Mordell}
L. J. Mordell,
\emph{On the rational solutions of the indeterminate equations of the
third and fourth degrees},
Proc. Camb. Phil. Soc. 21, 1922, pp. 179--192.

\bibitem{Ruck}
H. G. R\"uck,
\emph{A Note on Elliptic Curves Over Finite Fields},
Math. Comp. 49, 1987, pp. 301--304.

\bibitem{GrpSurvey}
M. Sala, and D. Taufer,
\emph{A survey on the group of points arising from elliptic curves with a Weierstrass model over a ring},
Int. J. Group Theory 12, 2023, pp. 177--196.

\bibitem{ELoops}
M. Sala, and D. Taufer,
\emph{Elliptic Loops},
J. Pure Appl. Algebra 227, 2023.

\bibitem{SatohAraki}
T. Satoh, and K. Araki,
\emph{Fermat quotients and the polynomial time discrete log algorithm for anomalous elliptic curves}, 
Comm. Math. Univ. Sancti Pauli 47, 1998, pp. 81--92.

\bibitem{Schoof}
R. Schoof,
\emph{Elliptic Curves over Finite Fields and the Computation of Square
Roots mod p},
Math. Comp. 44, 1985, pp. 483--494.

\bibitem{Semaev}
I. A. Semaev,
\emph{Evaluation of discrete logarithms in a group of p-torsion points
of an elliptic curve in characteristic p},
Math. Comp. 67, 1998, pp. 353--356. 

\bibitem{Shparlinski}
I. E. Shparlinski,
\emph{Pseudorandom number generators from elliptic curves},
Contemp. Math. 477, 2009, pp. 121--142.

\bibitem{Silverman}
J. H. Silverman,
\emph{The arithmetic of elliptic curves},
Springer-Verlag, 1986.

\bibitem{Smart}
N. Smart,
\emph{The discrete logarithm on elliptic curves of trace one},
J. Cryptology 12, 1999, pp. 193--196.

\bibitem{Voloch}
J. F. Voloch,
\emph{A note on elliptic curves over finite fields},
Bull. Soc. Math. France 116, 1988, pp. 455--458.

\bibitem{Washington}
L. C. Washington,
\emph{Elliptic curves, number theory and cryptography},
Chapman \& Hall / CRC, 2008.

\bibitem{Waterhouse}
W. C. Waterhouse,
\emph{Abelian varieties over finite fields},
Ann. Sci. \'Ecole Norm. Sup 2., 1969, pp. 521--560.

\bibitem{Wiles}
A. Wiles,
\emph{Modular elliptic curves and Fermat’s Last Theorem},
Ann. Math. 142, 1995, pp. 443--551.
\end{thebibliography}
%

\end{document}